\documentclass[10pt,twoside]{article}
\usepackage{lineno}
\usepackage{CJK}
\usepackage{graphicx}
\usepackage{amsfonts}
\usepackage{fancyhdr}
\usepackage{titlesec}
\usepackage{cite}
\usepackage{ifthen}
\usepackage{pifont}
\usepackage{stmaryrd}
\usepackage{setspace}
\usepackage{indentfirst}
\usepackage{amsmath,amssymb,amscd,bbm,amsthm,mathrsfs,dsfont}
\usepackage{bm}
\usepackage{authblk}
\input amssym.def

\newtheorem{conjecture}{Conjecture}[section]
\newtheorem{theorem}{Theorem}[section]

\newtheorem{lemma}{Lemma}[section]

\newtheorem*{pf}{Proof}
\theoremstyle{plain}

\newboolean{first}
\setboolean{first}{true}
\begin{document}
\title{\bf \Large Injective edge-coloring of graphs with small maximum degree
\footnotetext[1]{
\textbf{The first author's research is supported by NSFC (No. 12171436)}
}}
\author{Danjun Huang$^*$}
\author{Yuqian Guo}
\affil{\normalsize{Department of Mathematics, Zhejiang Normal University, Jinhua 321004, China}}
\date{}
\maketitle

\noindent
\textbf{Abstract}\quad An injective $k$-edge-coloring of a graph $G$ is a mapping $\phi$: $E(G)\rightarrow\{1,2,...,k\}$, such that $\phi(e)\ne\phi(e')$ if edges $e$ and $e'$ are at distance two, or are in a triangle. The smallest integer $k$ such that $G$ has an injective $k$-edge-coloring is called the injective chromatic index of $G$, denoted by $\chi_i'(G)$. In this paper, we prove that $\chi_i'(G)\le 7$ for every graph $G$ with $\Delta(G)\leq 4$ and mad$(G)<\frac{8}{3}$, where $\Delta(G)$ is the maximum degree of $G$.
\bigskip

\bigskip
\bigskip
\section{Introduction}
\par All graphs are simple, finite and undirected in this paper. For a graph $G$, we use $V(G)$, $E(G)$, $\Delta(G)$ and $\delta(G)$ to denote its vertex set, edge set, maximum degree and minimum degree, respectively. For a vertex $v\in V(G)$, let $N_G(v)$ be the set of neighbors of $v$ in $G$, and $d_G(v)=|N_G(v)|$ be the degree of $v$ in $G$. A vertex of degree $k$ (at least $k$, or at most $k$) is called a {\em$k$-vertex} ({\em$k^+$-vertex}, or {\em$k^-$-vertex}, respectively).
For a vertex $v\in V(G)$, we use $n_i^G(v)$ to denote the number of $i$-vertices adjacent to $v$ in $G$. For a graph $G$, the {\em maximum average degree} mad$(G)$ of $G$ is $\text{mad}(G)=\max\{\frac{2|E(H)|}{|V(H)|},H\subseteq G\}$.

\par An {\em injective $k$-edge-coloring} of a graph $G$ is a mapping $\phi$: $E(G)\rightarrow\{1,2,...,k\}$, such that $\phi(e)\ne\phi(e')$ if edges $e$ and $e'$ are at distance two, or are in a triangle. The smallest integer $k$ such that $G$ has an injective $k$-edge-coloring is called the {\em injective chromatic index} of $G$, denoted by $\chi_i'(G)$. The concept of injective edge-coloring was proposed in 2015 by Cardose et al. \cite{DJJC-2019} to slove the Packet Radio Network problem and they proved that it is NP-hard to compute the injective chromatic index of any graph. Moreover, Ferdjallah et al. \cite{BSA-2022} observed that the upper bound $\chi_i'(G) \leq 2(\Delta(G)-1)^2$ follows by Brooks' theorem. Meanwhile, they proposed the following conjecture.

\begin{conjecture}
For every subcubic graph $G$, $\chi_i'(G)\le 6$.
\end{conjecture}

In 2022, Miao et al. \cite{ZYG-2022} posed the following conjecture.

\begin{conjecture}
For every simple graph $G$ with maximum degree $\Delta$, $\chi_i'(G)\le \Delta(\Delta-1)$.
\end{conjecture}

\par Several authors have attacked this upper bound on the injective chromatic index for graphs with small maximum degree.  Towards Conjecture $1.1$, Kostochka et al. \cite{AAJ-2021} confirmed that $\chi_i'(G)\le 7$ and the upper bound $7$ can be improved to $6$ for subcubic planar graphs. Bu and Qi \cite{YC-2018} proved that $\chi_i'(G)\le 4$ for subcubic graphs with $\mathrm{mad}(G)<\frac{16}{7}$. Ferdjallah et al. \cite{BSA-2022} proved that $\chi_i'(G)\le 5$ for subcubic graphs with $\mathrm{mad}(G)<\frac{5}{2}$. Recently, Lu and Pan \cite{JX-2023} proved that $\chi_i'(G)\le 6$ for subcubic graphs with $\mathrm{mad}(G)<\frac{30}{11}$.

For the graphs $G$ with $\Delta(G)=4$,  Bu and Qi \cite{YC-2018} proved that $\mathrm{mad}(G)<m$ and $\chi_i'(G)\le k$ for some $(m,k)\in\{(\frac{10}{3},13),(\frac{18}{5},14),(\frac{15}{4},15)\}$. Lu and Pan \cite{JX-2023} proved that $\mathrm{mad}(G)<\frac{33}{10}$ and $\chi_i'(G)\le 12$. Miao et al. \cite{ZYG-2022} proved that $\mathrm{mad}(G)<m$ and $\chi_i'(G)\le k$ for some $(m,k)\in\{(\frac{14}{5},9),(3,10),(\frac{19}{6},11)\}$. In 2025, Fu and Lv \cite{JJ-2025} proved that $\mathrm{mad}(G)<m$ and $\chi_i'(G)\le k$ for some $(m,k)\in \{(\frac{5}{2},6),(\frac{13}{5},7),(\frac{36}{13},8)\}$.

\par Other related results on the injective edge coloring of graphs can be found in \cite{YPHJ-2024,JJN-2021,JYH-2023,JHY-2024}.

\par Motivated by the result of \cite{JJ-2025}, we consider the injective chromatic index of graphs with $\Delta(G)\le4$ and $\mathrm{mad}(G)<\frac{8}{3}$.

\begin{theorem}
Let $G$ be a graph with $\Delta(G)\leq 4$. If $\mathrm{mad}(G)<\frac{8}{3}$, then $\chi_i'(G)\le 7$.
\end{theorem}

\par Suppose that $G$ has a partial edge-coloring $\phi$ with the colors set $C$. For each edge $e'$ and $e$ in $G$, we say that edge $e'$ {\em sees} an edge $e$ if they are at distance two or are in a triangle. Moreover, $F(e)$ is the set of the edges that see $e$. A {\em$k_j$-vertex} is a $k$-vertex adjacent to $j$ $2$-vertices. A {\em $3_{1+}$-vertex} is a $3_1$-vertex adjacent to a $2_1$-vertex. A $3$-vertex $v$ is {\em poor} if $v$ is a $3_{1+}$-vertex or a $3_2$-vertex. The edge $e=uv$ is {\em light} if $d(u)=d(v)=2$.

\section{Proof of Theorem 1.1}
\par Let $G$ be a minimal counterexample to Theorem~1.1 with respect to $|V(G)| + |E(G)|$. Note that $G$ is connected with $\Delta(G) \leq 4$. Let $H$ be the graph obtained from $G$ by deleting all vertices of degree $1$. Since $H$ is the subgraph of $G$, mad$(H)\le$mad$(G)$. The following lemmas are proved by Fu and Lv in \rm{\cite{JJ-2025}}.

\begin{lemma}\rm{\cite{JJ-2025}} $(1)$ $\delta(H)\ge2$.
\par $(2)$ If $d_H(v)=2$, then $d_G(v)=2$.
\par $(3)$ If $v$ is a poor 3-vertex in $H$, then $d_G(v)=3$.
\end{lemma}

\begin{lemma}\rm{\cite{JJ-2025}}
 No $3$-vertex is adjacent to three poor 3-vertices in $H$.
\end{lemma}

\begin{lemma}\rm{\cite{JJ-2025}}
No $3_1$-vertex (or $3_2$-vertex) is adjacent to a poor 3-vertex in $H$.
\end{lemma}

\begin{lemma}\rm{\cite{JJ-2025}}
No $3_2$-vertex is adjacent to a $2_1$-vertex in $H$.
\end{lemma}

\textbf{Remark 2.1:}
Let $e=uv$ be the light edge in $G$ such that $N_G(u)=\{v,u'\}$ and $N_G(v)=\{u,v'\}$. Suppose that the subgraph $G'\subseteq G$ has an injective edge-coloring $\phi$ with $7$ colors and $e\in E(G')$. If $\phi(uu') \neq \phi(vv')$, then there are at most 3 edges incident to $v'$ that see $uv$ and at most 3 edges incident to $u'$ that see $uv$. So we may first erase the color of $uv$ and finally recolor it after arguing. In other words, we will omit the coloring of light edges in the following discussion.

\begin{lemma}
Let $C_3=vv_1v_2v$ be a triangle in $H$. If $d_H(v_1)=d_H(v_2)=2$, then $d_H(v)=4$, and $d_H(u)\ge 3$ for each $u\in N_H(v)\setminus\{v_1,v_2\}$.
\end{lemma}

\begin{pf}
\rm{Let $N_G(v)=\{v_1,v_2,u_1,\dots,u_{k-2}\}$. By Lemma $2.1$, $v_1v_2$ is a light edge in $G$.

\par Suppose that $d_H(v)\leq 3$. So $v$ is a $3_2$-vertex or a 2-vertex in $H$. By Lemma 2.1, $k=d_H(v)$.
By the minimality of $G$, $G'=G-v_1$ has an injective 7-edge-coloring $\phi$. Since $vv_1$ sees at most 3 edges
incident to $u_1$ if $u_1$ exists and $vv_1$ sees $vv_2$, we can extend $\phi$ to $G$, a contradiction.
So we prove that $d_H(v)=4$.

\par Now we need to show that $d_H(u_1)\geq 3$ and $d_H(u_2)\geq 3$. Suppose to the contrary that $d_H(u_1)=2$. By Lemma $2.1$, $d_G(u_1)=2$. By the minimality of $G$, $G'=G-v_1$ has an injective 7-edge-coloring $\phi$. Since $vv_1$ sees at most 4 edges
incident to $u_1$ or $u_2$, and $vv_1$ sees $vv_2$, we can extend $\phi$ to $G$, a contradiction. \qed
}
\end{pf}

\begin{lemma}
Any vertex $v$ in $G$ (and so in $H$) must have a $3^+$-neighbour.
\end{lemma}

\begin{pf}
\rm{Assume $G$ has a vertex $v$ with no $3^+$-neighbour. It has at most four neighbour, all $2^-$-vertices, say $v_i$ for $i\in\{1,2,3,4\}$. By Lemma 2.5 there are no edges between them. By the minimality of $G$, $G'=G-v$ has an injective 7-edge-coloring $\phi$. For each neighbour $v_i$ of $v$ there are at most 3 edges incident to $v$ that see $vv_i$ and at most 3 edges incident to the other neighbour $v_i'$ of $v_i$ that see $vv_i$. So we can extend $\phi$ to $vv_i$. Doing this for all neighbours $v_i$ of $v$, we can extend $\phi$ to $G$, a contradiction.\qed
}
\end{pf}

In the following discussion, if $v$ is the $k$-vertex in $H$, then we will use $v_1,v_2,\ldots,v_k$ to denote the $k$ neighbors of $v$ in $H$.

\begin{lemma}
 A $3$-vertex is adjacent to at most one $3_2$-vertex in $H$.
\end{lemma}
\begin{pf}
\rm{ Suppose that $v$ is adjacent to two $3_2$-vertices, say $v_1$  and $v_2$. By Lemma $2.1$, $d_G(v_1)=d_G(v_2)=3$. Let $N_G(v_i)=\{v,x_i,y_i\}$ for each $i\in\{1,2\}$. By Lemma $2.1$, $d_G(x_i)=d_G(y_i)=2$. By Lemma $2.5$, $x_1y_1\notin E(G)$ and $x_2y_2\notin E(G)$. Let $x_i'$ (or $y_i'$) be the second neighbor of $x_i$ (or $y_i$) other than $v_i$ for $i\in\{1,2\}$.

\textbf{Claim 2.1.} $d_G(v)=d_H(v)=3$.
\par\textbf{Proof of Claim 2.1.} Suppose that $v$ is adjacent to a $1$-vertex $u$ in $G$. By the minimality of $G$, $G'=G-u$ has an injective 7-edge-coloring $\phi$ with the color set $C$. Note that $F(vu)=\{v_ix_i,v_iy_i|\ i\in\{1,2\}\}\cup\{\phi(v_3w)|w\in N_G(v_3)\setminus\{v\}\}$ and $|F(vu)|\le7$. If $\phi$ cannot be extended to $vu$, then we may assume that $\phi(v_1x_1)=1$, $\phi(v_1y_1)=2$, $\phi(v_2x_2)=3$, $\phi(v_2y_2)=4$, and $\{\phi(v_3w)|w\in N_G(v_3)\setminus\{v\}\}=\{5,6,7\}$. Since $\{v_1x_1,v_1y_1\}\cup\{\phi(v_3w)|w\in N_G(v_3)\setminus\{v\}\}\subseteq F(vv_2)$, we have $S(vv_2)\subseteq\{3,4\}$, say $\phi(vv_2)=3$. If we can recolor $v_1x_1$ with a color $\alpha\in S(v_1x_1)\setminus\{1\}$, then we can extend $\phi$ to $G$  by letting $\phi(vu)=1$, a contradiction. Hence $F(v_1x_1)\cup\{1\}=C$. Since $|F(v_1x_1)|\le (d_G(x_1')-1)+(d_G(y_1)-1)+(d_G(v)-2)\le 3+1+2=6$, we have $3=\phi(vv_2)\notin F(v_1x_1)\setminus\{\phi(vv_2)\}$. Similarly, $F(v_1y_1)\cup\{2\}=C$ and $3=\phi(vv_2)\notin F(v_1y_1)\setminus\{\phi(vv_2)\}$. Note that $\phi(x_2x_2')\neq 4$. If $\phi(y_2y_2')\neq 4$, then we recolor or color $vv_2$, $v_1x_1$, $uv$ with $4$, $3$, $1$, respectively. The obtained coloring is the injective 7-edge-coloring of $G$, a contradiction. So $\phi(y_2y_2')=4$. Since $\{1,2\}\setminus\phi(x_2x_2')\neq\emptyset$, say $\phi(x_2x_2')\neq 1$, then we recolor or color $vv_2$, $v_1x_1$, $uv$ with $1$, $3$, $1$, respectively. The obtained coloring is the injective 7-edge-coloring of $G$, a contradiction.\qed

\par By the minimality of $G$, $G'=G-v_1$ has an injective 7-edge-coloring $\phi$. We erase the color of $v_2x_2$. Since $|F(vv_1)|\le (d_G(x_1)+d_G(y_1)-2)+(d_G(v_3)-1)+(d_G(v_2)-1)\le 6$, $|F(v_1x_1)|\le (d_G(x_1')-1)+(d_G(v)-1)+(d_G(y_1)-1)\le 6$, and $|F(v_1y_1)|\le (d_G(y_1')-1)+(d_G(v)-1)+(d_G(x_1)-1)\le 6$, we can color $vv_1,v_1x_1,v_1y_1$ independently. Now $|F(v_2x_2)|\le (d_G(v)-1)+(d_G(x_2')-1)+(d_G(y_2)-1)\le 6$. So we can extend $\phi$ to $G$, a contradiction.\qed
}
\end{pf}

\begin{lemma}
If a $3$-vertex $v$ is adjacent to a $3_2$-vertex in $H$, then $v$ is not adjacent to any $3_{1^+}$-vertex in $H$.
\end{lemma}
\begin{pf}
\rm{
\par Assume that $v$ is adjacent to a $3_2$-vertex $v_1$ and a $3_{1^+}$-vertex $v_2$. By Lemma 2.3, $v_1v_2\notin E(G)$. Let $N_H(v_i)=\{v,x_i,y_i\}$ for each $i\in\{1,2\}$ and $x_2$ is a $2_1$-vertex which is adjacent to a $2$-vertex $x_2'$ in $H$. By Lemma 2.1, $d_G(v_i)=3$ and $d_G(x_i)=d_G(y_1)=d_G(x_2')=2$ for each $i\in\{1,2\}$. Let $x_1'$ (or $y_1'$) be the second neighbor of $x_1$ (or $y_1$, respectively) other than $v_1$ .

\textbf{Claim 2.2.} $d_G(v)=d_H(v)=3$.
\par\textbf{Proof of Claim 2.2.} Suppose that $v$ is adjacent to a $1$-vertex $u$ in $G$. By the minimality of $G$, $G'=G-u$ has an injective 7-edge-coloring $\phi$ with the color set $C$. We erase the color of the light edge $x_2x_2'$.  Note that $F(vu)=\{v_1x_i,v_2y_i|i\in\{1,2\}\}\cup\{\phi(v_3w)|w\in N_G(v_3)\setminus\{v\}\}$ and $|F(vu)|\le7$. If $\phi$ cannot be extended to $vu$, then we may assume that $\phi(v_1x_1)=1$, $\phi(v_1y_1)=2$, $\phi(v_2x_2)=3$, $\phi(v_2y_2)=4$, and $\{\phi(v_3w)|w\in N_G(v_3)\setminus\{v\}\}=\{5,6,7\}$, which implies that $v_iv_j\notin E(G)$ for each $i,j\in\{1,2,3\}$. With a similar argument in Claim 2.1, we have that $S(vv_1)\subseteq\{1,2\}$ and $S(vv_2)\subseteq\{3,4\}$, say $\phi(vv_1)=1$ and $\phi(vv_2)=\beta\in \{3,4\}$. If we can recolor $v_1x_1$ with a color $\alpha \ne 1$, then we can extend $\phi$ to $G$  by letting $\phi(uv)=1$, a contradiction. Hence $F(v_1x_1)\cup\{1\}=C$. Since $|F(v_1x_1)|\le (d_G(x_1')-1)+(d_G(y_1)-1)+(d_G(v)-2)\le 3+1+2=6$, we have $\beta=\phi(vv_2)\notin F(v_1x_1)\setminus\{\phi(vv_2)\}$. Set $C_\phi(y_2)=\{\phi(y_2w)|w\in N_G(y_2)\setminus\{v_2\}\}$.

First suppose that $\beta=3$. If  $1\notin C_\phi(y_2)$, then we can recolor or color $vv_2$, $v_1x_1$, $uv$ with 1, 3, $1$, respectively, to extend $\phi$ to $G$, a contradiction. So $1\in C_\phi(y_2)$. Now $|F(v_2x_2)|\le (d_G(x_2')-1)+(d_G(y_2)-1)+1\le 1+3+1=5$. We recolor $v_2x_2$ with a color $\alpha\in S(v_2x_2)\setminus\{3\}$ and color $vu$ with $3$. The obtained coloring is the injective 7-edge-coloring of $G$, a contradiction. Now suppose that $\beta=4$. Note that $3\notin C_\phi(y_2)$. We recolor or color $vv_2$, $v_1x_1$, $uv$ with $3$, $4$, $1$, respectively. The obtained coloring is the injective 7-edge-coloring of $G$, a contradiction. \qed

\par By the minimality of $G$, $G'=G-x_2$ has an injective 7-edge-coloring $\phi$. Since $|F(v_2x_2)|\le (d_G(v)-1)+(d_G(y_2)-1)+(d_G(x_2')-1)\le 2+3+1=6$ and $x_2x_2'$ is light, we can extend $\phi$ to $G$ by Remark 2.1, a contradiction.\qed
}
\end{pf}

\begin{lemma}
Suppose that $v$ is a $4$-vertex adjacent to a $2_1$-vertex $v_1$ in $H$.
If each of $v_2$ and $v_3$ is a $2_0$-vertex or a $3_2$-vertex in $H$, then $d_G(v_4)=4$. In particular, $v_4$ is neither a $2$-vertex in $H$, nor a poor $3$-vertex in $H$.
\end{lemma}
\begin{pf}
\rm{Assume that $d_G(v_4)\le3$. Let $v_1$ be a $2_1$-vertex adjacent to a $2$-vertex $x_1$. By Lemma $2.1$, $d_G(v_1)=d_G(x_1)=2$, and $d_G(v_i)=d_H(v_i)$ for each $i\in\{2,3\}$.

Suppose that at least one of $v_2,v_3$ is a $2_0$-vertex, say $v_2$. By the minimality of $G$, $G'=G-v_1$ has an injective 7-edge-coloring $\phi$. Since $|F(vv_1)|\le (d_G(v_2)+d_G(v_3)+d_G(v_4)-3)+(d_G(x_1)-1)\le (2+3+3-3)+1=6$ and $v_1x_1$ is light, we can extend $\phi$ to $G$ by Remark 2.1, a contradiction.

\par Suppose that each $v_i$ is a $3_2$-vertex with $N_G(v_i)=\{v,y_{i1},y_{i2}\}$ for each $i\in\{2,3\}$. By Lemma $2.1$, $d_G(y_{i1})=d_G(y_{i2})=2$. Let $y_{i1}'$ $($or $y_{i2}'$$)$ be the second neighbor of $y_{i1}$ $($or $y_{i2}$, respectively$)$ other than $v_i$. By Lemma $2.5$, $y_{21}y_{22}\notin E(G)$. By the minimality of $G$, $G'=G-v_2$ has an injective 7-edge-coloring $\phi$. We erase the colors of $vv_3$ and $v_1x_1$. Since $|F(vv_2)|\le (d_G(v_3)+d_G(v_4)-2)+(d_G(y_{21})+d_G(y_{22})-2)=6$, $|F(v_2y_{21})|\le (d_G(y_{21}')-1)+(d_G(y_{22})-1)+(d_G(v)-2)\le 6$ and $|F(v_2y_{22})|\le (d_G(y_{22}')-1)+(d_G(y_{21})-1)+(d_G(v)-2)\le 6$, we color $vv_2$, $v_2y_{21}$, $v_2y_{22}$ independently. Now $|F(vv_3)|\le (d_G(v_2)+d_G(v_4)-2)+(d_G(y_{31})+d_G(y_{32})-2)\le 6$ and $v_1x_1$ is light, we can extend $\phi$ to $G$ by Remark 2.1, a contradiction.
\qed
}
\end{pf}

\begin{lemma}
 Suppose that $v$ is a $4$-vertex adjacent to two $2_1$-vertices in $H$, say $v_1$ and $v_2$.
Then each of $v_3$ and $v_4$ is neither a $2$-vertex in $H$, nor a poor 3-vertex in $H$.

\end{lemma}
\begin{pf}
\rm{Let $v_1$ and $v_2$ be two $2_1$-vertices adjacent to $v$. By Lemma $2.1$, $d_G(v_i)=2$ for each $i\in\{1,2\}$. Let $x_i$ be the second neighbor of $v_i$ other than $v$ for each $i\in\{1,2\}$. By Lemma $2.1$, $d_G(x_1)=d_G(x_2)=2$.

Assume that $v_3$  is a 2-vertex or a poor 3-vertex. By Lemma $2.1$, $d_G(v_3)=d_H(v_3)\le3$. By the minimality of $G$, $G'=G-v_1$ has an injective 7-edge-coloring $\phi$. We erase the color of $v_2x_2$. Since $|F(vv_1)|\le (d_G(v_3)+d_G(v_4)-2)+(d_G(x_1)-1)\le 6$ and both $v_1x_1$ and $v_2x_2$ are light, we can extend $\phi$ to $G$ by Remark 2.1, a contradiction.\qed
}
\end{pf}

\par Now we are ready to prove Theorem $1.1$. For each $v\in V(H)$, define the initial charge of $v$ as $w(v)=d_H(v)$. Moreover, we define appropriate discharging rules such that preserve the total charge sum of $H$. After the discharging procedure, we get the final charge $w^*(v)$ of $v$ is at least $\frac{8}{3}$. Thus $\frac{8}{3}|V(H)|\le \sum\limits_{v\in V(H)}{w^*(v)}=\sum\limits_{v\in V(H)}{w(v)}=\sum\limits_{v\in V(H)}{d_H(v)}< \frac{8}{3}|V(H)|$. This contradiction completes the proof of Theorem $1.1$.
\medskip
\par We now define the following discharging rules:

\medskip
\par\textbf{R1}~Each $3^+$-vertex sends $\frac{2}{3}$ to its adjacent $2_1$-vertex, and $\frac{1}{3}$ to its adjacent $2_0$-vertex.

\par\textbf{R2}~If $v$ is a non poor 3-vertex or a 4-vertex, then $v$ sends $\frac13$ to its adjacent $3_2$-vertex, and $\frac{1}{6}$ to its adjacent $3_{1^+}$-vertex.
\medskip

\par We need to check $\omega^*(v)$ for each vertex $v\in V(H)$. Let $d_H(v)=k$. By Lemma $2.1$, $k\ge 2$.

\textbf{Case 1}~$k=2$.

\par By Lemma $2.6$, $v$ is a $2_1$-vertex or a $2_0$-vertex. If $v$ is a $2_1$-vertex, then $\omega^*(v)=2+\frac{2}{3}=\frac{8}{3}$ by R1. If $v$ is a $2_0$-vertex, then $\omega^*(v)=2+\frac{1}{3}\times2=\frac{8}{3}$ by R1.

\textbf{Case 2}~$k=3$.

\par By Lemma $2.6$, $v$ is adjacent to at most two 2-vertices.

\par Suppose that $v$ is a $3_0$-vertex. By Lemma $2.7$, $v$ is adjacent to at most one $3_2$-vertex. If $v$ is adjacent to a $3_2$-vertex, then $v$ is not adjacent to any $3_{1^+}$-vertex by Lemma $2.8$. Hence, $\omega^*(v)\ge 3-\frac{1}{3}=\frac{8}{3}$ by R2. If $v$ is not adjacent to any $3_2$-vertex, then $v$ is adjacent to at most two $3_{1^+}$-vertices by Lemma $2.2$. Hence, $\omega^*(v)\ge 3-2\times\frac{1}{6}=\frac{8}{3}$ by R2.

\par Suppose that $v$ is a $3_1$-vertex. By Lemma $2.3$, $v$ is not adjacent to any poor 3-vertex. If $v$ is a $3_{1^+}$-vertex, then $\omega^*(v)\ge 3-\frac{2}{3}+2\times\frac{1}{6}=\frac{8}{3}$ by R1 and R2. If $v$ is not a $3_{1^+}$-vertex,  then $\omega^*(v)\ge 3-\frac{1}{3}=\frac{8}{3}$ by R1 and R2.

\par Suppose that $v$ is a $3_2$-vertex. By Lemma $2.3$ and Lemma 2.4, $v$ is adjacent to two $2_0$-vertices and is not adjacent to any poor 3-vertex.  Hence, $\omega^*(v)\ge 3-2\times\frac{1}{3}+\frac{1}{3}=\frac{8}{3}$ by R1 and R2.

\textbf{Case 3}~$k=4$.

\par By Lemma $2.6$, $v$ is adjacent to at most three $2$-vertices. If $v$ is not adjacent to any $2_1$-vertex, then $\omega^*(v)\ge 4-4\times\frac{1}{3}=\frac{8}{3}$ by R1 and R2. If $v$ is adjacent to two $2_1$-vertices, then $v$ is not adjacent to any poor 3-vertex,  and $v$ is adjacent to two $3^+$-vertices by Lemma $2.10$. Hence, $\omega^*(v)=4-2\times\frac{2}{3}=\frac{8}{3}$ by R1 and R2. So suppose that $v$ is adjacent to exactly one $2_1$-vertex. Let $s(v)$ be the sum of the number of $2_0$-vertices and $3_2$-vertices adjacent to $v$. By Lemma 2.9, $s(v)\le 2$.
If $s(v)\le1$, then $\omega^*(v)=4-\frac{2}{3}-\frac{1}{3}\times s(v)-\frac{1}{6}\times(3-s(v))=\frac{17}{6}-\frac{1}{6}s(v)\ge \frac{8}{3}$ by R1 and R2. If $s(v)=2$, then $v$ is not adjacent to any $3_{1^+}$-vertex by Lemma $2.9$. Hence, $\omega^*(v)=4-\frac{2}{3}-2\times\frac{1}{3}=\frac{8}{3}$ by R1 and R2.

\end{document}